\documentclass[11pt, reqno]{amsart}
\usepackage{amsmath,amsthm,amssymb}
\theoremstyle{plain}
\newtheorem{thm}{Theorem}[section]
\newtheorem{prop}[thm]{Proposition}
\newtheorem{lem}[thm]{Lemma}
\newtheorem{cor}[thm]{Corollary}

\theoremstyle{definition}

\newtheorem{rem}[thm]{Remark}
\newtheorem{defn}[thm]{Definition}
\newtheorem{eg}[thm]{Example}
\newtheorem{subtitle}[thm]{}
\newtheorem{ex}{Exercise}[section]
\numberwithin{equation}{section}

\def\d{\delta}

\def\e{\epsilon}

\def\l{\lambda}

\def\n{\,\vert\,}
\def\o{\theta}

\def\W{\Omega}

\def\bs{\bigskip}
\def\ms{\medskip}
\def\ss{\smallskip}

\def\ni{\noindent}
\def\ti{\tilde}
\def\p{\partial}

\def\I{{\rm I\/}}
\def\II{{\rm II\/}}
\def\diag{{\rm diag}}
\def\ad{{\rm ad}}
\def\Ad{{\rm Ad}}

\def\R{\mathbb{R} }
\def\C{\mathbb{C}}

\newcommand{\beg}{\begin{eg}}
\newcommand{\eeg}{\end{eg}}
\newcommand{\bthm}{\begin{thm}}
\newcommand{\ethm}{\end{thm}}
\newcommand{\bprop}{\begin{prop}}
\newcommand{\eprop}{\end{prop}}
\newcommand{\bcor}{\begin{cor}}
\newcommand{\ecor}{\end{cor}}
\newcommand{\blem}{\begin{lem}}
\newcommand{\elem}{\end{lem}}
\newcommand{\bca}{\begin{cases}}
\newcommand{\eca}{\end{cases}}
\newcommand{\brem}{\begin{rem}}
\newcommand{\erem}{\end{rem}}
\newcommand{\bpm}{\begin{pmatrix}}
\newcommand{\epm}{\end{pmatrix}}
\newcommand{\bbm}{\begin{bmatrix}}
\newcommand{\ebm}{\end{bmatrix}}
\newcommand{\bvm}{\begin{vmatrix}}
\newcommand{\evm}{\end{vmatrix}}
\newcommand{\bdefn}{\begin{defn}}
\newcommand{\edefn}{\end{defn}}
\newcommand{\bsub}{\begin{subtitle}}
\newcommand{\esub}{\end{subtitle}}
\newcommand{\bex}{\begin{ex}}
\newcommand{\eex}{\end{ex}}
\newcommand{\ben}{\begin{enumerate}}
\newcommand{\een}{\end{enumerate}}

\newcommand{\balign}{\begin{align}}
\newcommand{\ealign}{\end{align}}
\newcommand{\baligns}{\begin{align*}}
\newcommand{\ealigns}{\end{align*}}
\newcommand{\beq}{\begin{equation}}
\newcommand{\eeq}{\end{equation}}

\def\calA{{\mathcal A}}
\def\calB{{\mathcal B}}
\def\calG{{\mathcal G}}
\def\calJ{{\mathcal J}}
\def\calK{{\mathcal K}}
\def\calL{{\mathcal L}}

\def\calP{{\mathcal P}}
\def\calS{{\mathcal S}}
\def\calU{{\mathcal U}}
\def\calQ{{\mathcal Q}}
\def\calO{{\mathcal O}}

\def\({\left(}
\def\){\right)}
\def\[{\left[}
\def\]{\right]}

\def\rd{{\rm d\/}}

\def\onn{\frac{O(2n)}{ O(n)\times O(n)}}

\begin{document}

\title[Soliton Hierarchies from Involutions]
{ Soliton Hierarchies Constructed from Involutions}
\author{Chuu-Lian Terng$^\dag$}\thanks{$^\dag$Research supported
in  part by NSF Grant DMS-0707132}
\address{ Department of Mathematics\\ University of California, Irvine, CA 92697-3875}
  \email{cterng@math.uci.edu}

\begin{abstract}

We introduce two families of soliton hierarchies: the twisted hierarchies associated to symmetric spaces.  The Lax pairs of these two hierarchies are Laurent polynomials in the spectral variable.  Our constructions gives a hierarchy of commuting flows for the generalized sine-Gordon equation (GSGE), which is the Gauss-Codazzi equation for $n$-dimensional submanifolds in Euclidean $2n-1$ space with constant sectional curvature $-1$. In fact, the GSGE is the first order system associated to a twisted Grassmannian system.  We also study symmetries for the GSGE.

\vskip 4.5mm

\noindent {\bf 2000 Mathematics Subject Classification:} 37K05, 37K25, 53B25.

\noindent {\bf Keywords and Phrases:} Soliton hierarchy, Involutions, Submanifolds.
\end{abstract}

\vskip 12mm

\maketitle

\section{Introduction} \label{section 1}

Most soliton hierarchies can be constructed from splittings of Lie algebras $\calL$ as positive $\calL_+$ and negative  $\calL_-$ subalgebras (cf. \cite{ZS72, AKNS74, Adl79, DriSok84, TerUhl00a}).  For example, the AKNS $2\times 2$,  the Korteweg-de-Vris (KdV), the non-linear Schr\"odinger (NLS), the $3$-wave, and the Gelfand-Dikki hierarchies are constructed in this way.  When the infinite dimension Lie algebra is a subalgebra of the Lie algebra $\calL(\calG)$ of loops in a complex, simple Lie algebra $\calG$, many properties of soliton equations can be explained in a unified way including: Lax pairs, commuting flows, bi-Hamiltonian structures, B\"acklund transformations, scattering and inverse scattering,  tau functions, and the Virasoro actions (cf. \cite{Wil91, TerUhl98, TerUhl07a, TerUhl09a, TerUhl09b}).

  If we use the standard splitting of $\calL(\calG)$ as non-negative and negative Fourier series, then we get the  $G$-hierarchy \cite{Wil91}.  If $U$ is the real form of $G$ defined by an involution $\tau$, then $\tau$ induces an involution on $\calL(\calG)$ and its fixed point set $\calL^\tau(\calG)$ in $\calL(\calG)$ leaves the positive and negative loop algebras invariant. We call the resulting hierarchy the $U$-hierarchy. For example,  the $SU(2)$-hierarchy contains NLS. 
   Assume further that there is a second involution $\sigma$ commuting with $\tau$.  Let $K$ denote the fixed point set of $\sigma$ in $U$. Then $\frac{U}{K}$ is a symmetric space.  There is also an induced involution given by $\sigma$ on $\calL(\calG)$ and the subalgebra $\calL^{\tau,\sigma}(\calG)$ that is fixed by the involutions induced by $\tau$ and $\sigma$ gives the $\frac{U}{K}$-hierarchy.  For example, the $SU(2)/SO(2)$-hierarchy contains the modified KdV (cf. \cite{TerUhl00a}).

If the soliton hierarchy is constructed from a subalgebra of $\calL(\calG)$, then we can put all degree $1$ flows together to form a first order system of partial differential equations. Here a soliton flow is said to have degree one  if its Lax pair is a degree one polynomial in the spectral (or loop) variable. We call the first order system constructed from the $U$- and $\frac{U}{K}$- hierarchies the $U$- and the $\frac{U}{K}$- systems respectively (\cite {Ter97}). Techniques in soliton theory can be applied to these first order systems.     Many of the $\frac{U}{K}$- systems are the Gauss-Codazzi equations for submanifolds in space forms or symmetric spaces with special geometric properties.  For example, the Gauss-Codazzi equations for isothermic surfaces in $\R^3$ is the $\frac{O(4,1)}{O(3)\times O(1,1)}$-system (\cite{CGS95}), for isothermic surfaces in $\R^n$ is the $\frac{O(n+1,1)}{O(n)\times O(1,1)}$-system (\cite{BDPT02, Bur06}),  for $n$-dimensional flat submanifolds in $\R^{2n}$ with flat and non-degenerate normal bundle is the $\frac{O(2n)}{O(n)\times O(n)}$-system (\cite{Ter97}), for flat Lagrangian submanifolds in $\C^n=\R^{2n}$ with non-degenerate normal bundle is the $\frac{U(n)\ltimes \C^n}{O(n)\ltimes \R^n}$-system (\cite{TerWan06}),  for flat Lagrangian submanifolds in $\C P^{n-1}$ is the $\frac{U(n)}{O(n)}$-system, and the equation for conformally flat $n$-submanifolds in $S^{2n-2}$ with flat and non-degenerate normal bundle is the $\frac{O(2n-1, 1)}{O(n)\times O(n-1,1)}$-system (\cite{DoTe07}).

The (twisted) $\frac{U}{K}$-system is the equation for regular curved flats in $\frac{U}{K}$  studied by Ferus and Pedit in \cite{FP96a} written in a good choice of coordinate system.  Finite type curve flats were studied by Ferus and Pedit in \cite{FP96a, FP96b} and by Brander in \cite{Bra07a, Bra07b}.  

It is known that the $\frac{U}{K}$-system as an exterior differential system is involutive, hence by the Cartan-K\"ahler Theorem the Cauchy problem has a unique local real analytic solution for real analytic initial data on a non-characteristic line (cf. \cite{TerWan05}).  We can also use the inverse scattering method to prove that the Cauchy problem has unique global smooth solution for smooth rapidly decaying small initial data on a non-characteristic line \cite{Ter97}. As a consequence, given any smooth initial data defined on an interval $\I$ of a non-characteristic line, there is a subinterval $\I_0\subset \I$ and a solution of the $\frac{U}{K}$-system that agrees with the initial data on $\I_0$. 

The group of rational maps $f:S^2\to U_\C$ satisfying the $\frac{U}{K}$-reality condition and $f(\infty)=\I$ acts on the space of solutions of the $\frac{U}{K}$-system, and solutions in the orbit at the vacuum solution are given explicitly in terms of exponential functions (\cite{TerUhl00a}). 

Since the Gauss-Codazzi equation for surfaces in $\R^3$ with constant curvature $-1$ is the sine-Gordon equation,  the Gauss-Codazzi equation for $n$-dimensional submanifolds in $\R^{2n-1}$ with constant sectional curvature $-1$ is called the {\it generalized sine-Gordon equation\/} (GSGE).  The GSGE has B\"acklund transformations, Permutability formula (\cite{TenTer80, Ter80}), a Lax pair, and the inverse scattering theory has been carried out in \cite{AbBeTe86, BeTe88}.   But it is not clear whether there exists a soliton hierarchy constructed from a splitting of some subalgebra of a loop algebra so that the first order system associated to this hierarchy is the GSGE.  If we can find such splitting, then techniques from soliton theory will give us infinitely many commuting flows and a loop group action on the space of solutions of the GSGE. The Lax pairs of the $U$- and $\frac{U}{K}$-systems are degree $1$ polynomials in $\l$. But  the Lax pair for the GSGE is a degree $1$ Laurent polynomial in $\l$.  So the GSGE cannot be a $U$- or $\frac{U}{K}$-system.  This makes the problem of finding a splitting that gives GSGE more interesting because as a bi-product this should give a way to construct new soliton hierarchies.  In fact, this is what we have done in this paper.  If $\frac{U}{K}$ is a non-compact symmetric space given by $\tau, \sigma$,  then we can use $\sigma$ to construct a splitting of the loop algebra $\calL^\tau(\calG)$  and use it to build a twisted $U$-hierarchy.  If there is a third involution $\sigma_2$ that commutes with $\tau, \sigma$ and these three involutions satisfy certain conditions, then we can use $\sigma_2$ to construct a splitting of the loop algebra $\calL^{\tau,\sigma}(\calG)$ and use it to construct a twisted $\frac{U}{K}$-hierarchy.  Lax pairs of flows in these two  twisted hierarchies are Laurent polynomials in the spectral variable.

This paper is organized as follows: In section 2, we review the method of constructing soliton hierarchies from splittings of Lie algebras and use these splittings to explain some soliton properties.  In section 3, we first review the construction of the $G$-, $U$-, and $\frac{U}{K}$-hierarchy of soliton equations, then give constructions of two new soliton hierarchies,  the twisted $U$- and twisted $\frac{U}{K}$- hierarchies. In the final section, we review the construction of the $U$-system and the $\frac{U}{K}$-system and use the same method to construct the twisted $U$- and twisted $\frac{U}{K}$-system. We also show that the GSGE is the first order system associated to a twisted $\frac{O(n,n)}{O(n)\times O(n)}$-hierarchy.  As a bi-product, we obtain a loop group action and a hierarchy of commuting flows on the space of solutions of GSGE.

\section{Splittings and soliton hierarchies}

We review the method for constructing soliton hierarchies from splittings of Lie algebras and vacuum sequences and give several examples ((cf. \cite{AKNS74, DriSok84, Wil91, TerUhl00a, TerUhl07a}).  

\ss\ni {\bf Lie algebra splitting and vacuum sequence} 

Let $G$ be a complex simple Lie group,  $L(G)$ the group of smooth maps $f:S^1\to G$, $\calG$ the Lie algebra of $G$, and  $\calL(\calG)$ the Lie algebra of $L(G)$.  Suppose $\calL$ is a Lie subalgebra of $\calL(\calG)$.  A {\it splitting of $\calL$} is a pair of  Lie subalgebras $(\calL_+, \calL_-)$ of $\calL$ such that $\calL=\calL_+\oplus \calL_-$ as direct sum of linear subspaces and the corresponding subgroups $L_+, L_-$ intersect at the identity (i.e., $L_+\cap L_-=\{e\}$). 

A commuting sequence $\calJ=\{J_j\n j\geqslant 1\}$ in $\calL_+$ is called a {\it vacuum sequence\/} if they are linearly independent and $J_j$ is an analytic function of $J_1$ for all $j\geqslant 1$.

\ss\ni {\bf Construction of flows}

A soliton hierarchy can be constructed from the splitting $(\calL_+, \calL_-)$ of $\calL$ and a vacuum sequence $\calJ$ (cf. \cite{TerUhl98, TerUhl07a}) as follows:  Set
\beq\label{ak}
Y= \{\pi_+(gJ_1g^{-1})\n g\in L_-\},
\eeq
where $\pi_\pm:\calL\to \calL_\pm$ are the projections with respect to $\calL= \calL_+ \oplus\calL_-$. Note that $J_1\in Y$. 

 Assume that given smooth $P:\R\to Y$, there exists $M:\R\to Y$ such that 
 $$\p_x+ P = M(\p_x+ J_1)M^{-1},$$
 or equivalently, 
 $$P= MJ_1M^{-1} - M_xM^{-1}.$$
 Such $M$ is called a {\it Baker function for $P$}.  
 Since $P\in \calL_+$ and $M_xM^{-1}\in \calL_-$, we have
 \beq\label{bx}
 P= \pi_+(MJ_1M^{-1}), \quad M_xM^{-1}= \pi_-(MJ_1M^{-1}).
 \eeq
 Then  $Q_j(P)= MJ_jM^{-1}\in \calL$ satisfies the following properties:
\beq\label{bt} 
\bca [\p_x+P, Q_j(P)]=0,\\
Q_j(P)\,\, {\rm  is\,  conjugate\, to\, \, }J_j,\\
Q_j(J_1)= J_j.
\eca
\eeq 
In fact, $Q_j(P)$ is determined by these properties.  

The flow generated by $J_j$ is the equation on $C^\infty(\R, Y)$:
\beq\label{aa}
P_{t_j}= [\p_{x} +P, \, \pi_+(Q_j(P))].
\eeq
These flows commute and the collection of these commuting flows is {\it the soliton hierarchy associated to the splitting $(\calL_+, \calL_-)$ and the vacuum sequence $\calJ$\/}.  If $Q_j(P)$ depends on $P$ and derivatives of $P$, then these flows are PDEs.

\ss
\ni{\bf Lax pair and frame\/}

The following statements are equivalent for smooth maps $P:\R^2\to Y$:
\ben 
\item $P$ is a solution of the flow equation \eqref{aa} generated by $J_j$,
\item $[\p_x+P, \, \p_{t_j} + \pi_+(Q_j(P))]=0$,
\item $\o= P\rd x+ \pi_+(Q_j(P))\rd t_j$ is a flat connection $1$-form,
\item 
\beq\label{ca}
\bca E^{-1} E_x= P, \\
 E^{-1}E_{t_j}= \pi_+(Q_j(P)),
 \eca
 \eeq
 is solvable for any initial data $E(0,0)=c_0\in L_+$.
 \een
 We call $\o$ a {\it Lax pair\/}, and call the solution $E$ of \eqref{ca} with initial data $c_0=$ the identity in $L_+$ the {\it frame\/} of the solution $P$ of the flow generated by $J_j$.

Note that the constant $P=J_1$ is a solution of the $j$-th flow and its frame is $V(x,t_j)=\exp(xJ_1 + t_j J_j)$, which are called the {\it vacuum solution\/} and the {\it vacuum frame\/} respectively. 

\bthm[Local Factorization Theorem]\label{bv} \cite{PreSeg86, TerUhl09a}\par

Suppose $L$ is a closed subgroup of the group of Sobolev $H^1$- loops in a finite dimensional Lie group $G$, and $(\calL_+, \calL_-)$ is a splitting of the Lie algebra $\calL$.  Let $\calO$ be an open subset in $\R^N$, and $g:\calO\to L$ a map such that $(x, \l)\mapsto g(x)(\l)$ is smooth.  If $p_0\in \calO$ and $g(p_0)= k_+k_-= h_-h_+$ with $k_\pm, h_\pm\in L_\pm$, then there exist an open subset $\calO_0\subset \calO$ containing $p_0$ and unique $f_\pm, g_\pm : \calO_0\to L_\pm$ such that $g= g_+g_-= f_-f_+$ on $\calO_0$ and $g_\pm(p_0)= k_\pm$, $f_\pm(p_0)= h_\pm$.
\ethm

\ss\ni{\bf Formal inverse scattering} \cite{TerUhl07a}

Given $f\in L_-$, we can use The Local Factorization Theorem \ref{bv} to construct a local solution $P_f$ of the hierarchy as follows:
Let
$$V(x,t_j)=\exp(xJ_1+ t_jJ_j)$$
denote the vacuum frame.
By Theorem \ref{bv} there is an open subset $\calO$ of the origin in $\R^2$ such that we can factor $fV(x,t_j)$ as
\beq\label{ab}
fV(x,t_j)=E(x,t_j) M(x,t_j)
\eeq
with $E(x,t_j)\in L_+$ and $M(x,t_j)\in L_-$ for all $(x, t_j)\in \calO$.
Then \beq\label{ac}
P_f= \pi_+(MJ_1M^{-1})
\eeq
 is a solution of the $j$-th flow \eqref{aa}.
 
 \ss
\ni {\bf Dressing action and B\"acklund transformations} (\cite{TerUhl00a})

Given a solution $P$ of the flow generated by $J_j$ and $f\in L_-$, let $E$ denote the frame of $P$.  Then by Theorem \ref{bv} there is an open subset 
$\calO$ of the origin in $\R^2$ such that we can factor
$$fE(x,t_j)= \ti E(x,t_j)\ti f(x,t_j)$$
with $\ti E(x, t_j)\in L_+$ and $\ti f(x,t_j)\in L_-$ pointwise for each $(x,t_j)\in \calO$. Moreover, 
\ben
\item $\ti P= \ti E^{-1} \ti E_x$ is a solution of the flow generated by $J_j$ and $\ti E$ is the frame of $\ti P$,
\item $f\ast P=\ti P$ defines an action of $L_-$ on the space of local solutions of the flow generated by $J_j$, 
\item if $L$ is a subgroup of the group of smoth loops in $GL(n)$ and $f\in L_+$ is rational, then
\ben
\item $P\mapsto f\ast P$ can be computed explicitly using the frame of $P$ and poles and residues of $f$,
\item if  $f$ is a {\it simple element\/}, i.e., $f$ can not be factored as the product of two rational elements in $L_-$, then $f\ast P$ can be also computed by solving a system of compatible ODEs, which are usually called {\it B\"acklund or Darboux transformations\/},
\item $f\ast J_1$ is given explicitly in terms of exponential functions and is a pure soliton soluton.
\een
\een

\section{Examples}

\ni {\bf The $G$-hierarchy }\,{\rm (\cite{AKNS74, Wil91})}\par

\ss
Let $L(G)$ denote  the group of smooth loops $f:S^1\to G$, $L_+(G)$ is the subgroup of $f\in L(G)$ such that $f$ is the boundary value of a holomorphic map defined on $|\l |<1$, and $L_-(G)$ is the subgroup of $f\in L(G)$ such that $f$ is the boundary value of a holomorphic map $h$ defined on $1< |\l |\leqslant \infty$ and $h(\infty)=e$ the identity in $G$.  It follows from the Liouville Theorem that $L_+(G)\cap L_-(G)=\{e\}$. 
Let $\calL(\calG), \calL_\pm(\calG)$ be the corresponding Lie algebras. It is clear that $\calL(\calG)=\calL_+(\calG)\oplus \calL_-(\calG)$.  So
$(\calL_+(\calG), \calL_-(\calG))$ is a splitting of $\calL(\calG)$.

%:
Let $\calA$ be a Cartan subalgebra of $\calG$, and $\{a_1, \ldots, a_n\}$ a basis of $\calA$ such that $a_1$ is regular, i.e., the centralizer of $a_1$, $\{x\in \calG\n [x,a_1]=0\}$ is equal to $\calA$.  Then $\calJ=\{J_{i,j}=a_i \l^j\n 1\leqslant i\leqslant n, \, j\geqslant 1\}$ is a vacuum sequence.  The soliton hierarchy constructed from the splitting $(\calL_+(\calG), \calL_-(\calG))$ and the vacuum sequence $\calJ$ is called the $G$-hierarchy.

\ss\ni {\bf The $U$-hierarchy} \,{\rm (\cite{TerUhl00a})}

Let $\tau$ be an involution of $G$ such that the differential at the identity $e\in G$ (still denoted by $\tau$) is conjugate linear.  The fixed point set $\calU$ of $\tau$ in $\calG$ is called a {\it real form\/} of $\calG$.
Let $\calA$ be a maximal abelian subalgebra in $\calU$, and $\{a_1, \ldots, a_n\}$ a basis of $\calA$ such that the centralizer of $a_1$ in $\calU$ is $\calA$. Let $\calL(\calG)$ and $\calL_\pm(\calG)$ be as in the $G$-hierarchy.  Set
\begin{align*}
\calL^{\tau}(\calG)&=\{\xi\in \calL(\calG)\n \tau(\xi(\bar\l))= \xi(\l)\}, \\
\calL^{\tau}_\pm(\calG)&=\calL_\pm(\calG)\cap \calL^\tau(\calG),\\
\calJ&= \{ J_{i,j}=a_i \l^j\n 1\leqslant i\leqslant n, \, j\geqslant 1\}.
\end{align*}
Then $(\calL^\tau_+(\calG), \calL^\tau_-(\calG))$ is a splitting of $\calL^\tau(\calG)$, and $\calJ$ is a vacuum sequence.  The hierarchy generated by this splitting and vacuum sequence is called the $U$-hierarchy.  For example, for $G=SL(2,\C)$ and $\tau(x)= (\bar x^t)^{-1}$, we have $U=SU(2)$.  Choose $a_1=\diag(i, -i)$.  The space $Y$ defined by \eqref{ak} is
$$Y=\left\{a\l + \bpm 0& q\\ -\bar q& 0\epm\,\bigg|\, q\in \C\right\},$$
and the flows generated by $a_1\l^2$ is the NLS, 
$ q_t= \frac{i}{2}(q_{xx}+ 2|q|^2 q)$.

   For general $U$, the space $Y=\{ a_1\l+\xi \n \xi\in\calA^\perp\}$, where $\calA^\perp=\{\xi\in \calU\n (\xi, \calA)=0\}$ and $(\, , )$ is a fixed non-degenerate bilinear form on $\calU$.  The  first flow generated by $a_i\l$ is 
   \beq\label{aj}
   u_{t_i}= \ad(a_i)\ad(a_1)^{-1}(u_x) + [u, \ad(a_i)\ad(a_1)^{-1}(u)],
   \eeq
   where $\ad(b)(x)=[b,x]$ (since $a_1$ is regular, $\ad(a_1)$ is a linear isomorphism on $\calA^\perp$).
  
\ss\ni {\bf The $\frac{U}{K}$-hierarchy} \, {\rm (\cite{TerUhl00a})}

Let $\tau, \sigma$ be two commuting involutions of $G$ such that the differentials $d\tau_e$ is conjugate linear and $d\sigma_e$ is complex linear on $\calG$. To simplify the notations, we will use $\tau, \sigma$ to denote the induced involutions $d\tau_e, d\sigma_e$ on $\calG$.  Let $U$ be the fixed point set of $\tau$ in $G$, and $K$ the fixed point set of $\sigma$ in $U$. The quotient space $\frac{U}{K}$ is a symmetric space.  Let $\calP$ denote the $-1$ eigenspace of $\sigma$ on $\calU$.  Then $\calU=\calK+ \calP$ is a {\it Cartan decomposition\/} of $\frac{U}{K}$.  Let $\calA$ be a maximal abelian subalgebra in $\calP$.  The dimension of $\calA$ is the rank of the symmetric space.   Let $\{a_1, \ldots, a_n\}$ be a basis of $\calA$ such that $a_1$ is regular, i.e.,  the Ad$(K)$ orbit at $a_1$ in $\calP$ is a maximal orbit.  Set
\begin{align*}
\calL^{\tau, \sigma}(\calG)&=\{\xi\in \calL(\calG)\n \tau(\xi(\bar\l))=\xi(\l), \, \sigma(\xi(-\l))=\xi(\l)\},\\
\calL^{\tau, \sigma}_\pm(\calG)&= \calL^{\tau,\sigma}(\calG)\cap \calL_\pm(\calG),\\
\calJ&=\{J_{i, 2k+1}= a_i \l^{2k+1}\n 1\leqslant i\leqslant n, \, k\geqslant 0\}.
\end{align*}
Then $(\calL^{\tau,\sigma}_+(\calG), \calL^{\tau, \sigma}_-(\calG))$ is a splitting of $\calL^{\tau, \sigma}(\calG)$ and $\calJ$ is a vacuum sequence.  For example, let $G= SL(2,\C)$, $\tau(g)= (\bar g^t)^{-1}$, $\sigma(g)= (g^t)^{-1}$, and $a_1=\diag(i, -i)$.  The symmetric space is $\frac{SU(2)}{SO(2)}\simeq \C P^1$, the space
$$Y=\left\{a_1\l + \bpm 0 &q\\ -q&0\epm \,\bigg|\, q: \R\to \R\right\}.$$
The flow generated by $a_1\l^3$ in the $\frac{SU(2)}{SO(2)}$-hierarchy is the mKdV,
$q_t= \frac{1}{3} (q_{xxx} + 6 q^2q_x)$.

\ss\ni {\bf The $U$-hierarchy twisted by $\sigma$}

Let $\frac{U}{K}$ be the symmetric space defined by involutions $\tau$ and $\sigma$ of $G$.  Suppose
\ben
\item[(a)] the real form $U$ is non-compact and there is a Borel subgroup $B$ such that
$U=KB, \quad K\cap B=\{e\}, \quad \calU= \calK + \calB$,
\item[(b)] there is a maximal abelian subalgebra $\calA$ in $\calU$ such that $\sigma(\calA)=\calA$.
\een
Below we associate to such a symmetric space a splitting and a vacuum sequence. 

\ss
Let $\ti L^\tau(G)$ be the group of smooth maps $f:S^1_{\e^{-1}}\to G$ satisfying $\tau(f(\bar\l))= f(\l)$, where $S^1_{\e^{-1}}=\l\in \C\n |\l |=\e^{-1}\}$. Let $\ti L^\tau_+(G)$ be the subgroup of $f\in \ti L^\tau(G)$ such that there exists a holomorphic map $h$ on $\e < |\l| <\e^{-1}$ with boundary values $f$ on $S^1_{\e^{-1}}$ and $\hat f(\l)= \sigma(f(\l^{-1}))$ on $S^1_\e$, and $\ti L^\tau_-(G)$ the subgroup of $f\in \ti L^\tau(G)$ such that there exists a holomorphic map $h$ defined on $\infty\geqslant |\l| >\e^{-1}$ with boundary value $f$ on $S^1_{\e^{-1}}$ and $h(\infty)\in B$. It follows from the Liouville Theorem that $\ti L_+^\tau(G)\cap \ti L^\tau_-(G)= \{e\}$. 
 The corresponding Lie subalgebras are:
\begin{align*}
&\ti\calL^\tau(\calG)= \{\xi(\l)= \sum_{j\leqslant n_0}\xi_j\l^j\n \xi_j\in \calU\},\\
&\ti\calL_+^\tau(\calG)=\{\xi\in \ti\calL^\tau(\calG)\n \sigma(\xi(\l^{-1}))= \xi(\l)\},\\
&\ti\calL_-^\tau(\calG)=\{\xi\in \ti\calL^{\tau}(\calG)\n \xi(\l)=\sum_{j\leqslant 0} \xi_j \l^j, \,\, \xi_0\in \calB\}.\\
\end{align*}
 So $(\ti\calL_+^\tau(\calG), \ti \calL^\tau_-(\calG))$ is a splitting of $\ti\calL^\tau(\calG)$. 
 Let $\pi_\pm:\ti\calL^{\tau}(\calG)\to \ti\calL_\pm^{\tau}(\calG)$ denote the projection with respect to $\ti\calL^{\tau}(\calG)=\ti \calL_+^{\tau}(\calG)\oplus \ti \calL_-^{\tau}(\calG)$.  In fact, for $A(\l)=\sum_{j\leqslant n_0} A_j \l^j$ in $\ti\calL^\tau(\calG)$, we have
\begin{align*}
&\pi_+(A)= (A_0)_{\calK} + \sum_{j>0} (A_j \l^j+\sigma(A_j)\l^{-j}),\\
& \pi_-(A)= (A_0)_{\calB} + \sum_{j>0} (A_{-j} -\sigma(A_j))\l^{-j},
\end{align*}
where $(A_0)_{\calK}$ and $(A_0)_{\calB}$ are the projections of $\calU$ onto $\calK$ and $\calB$ respectively.

  Let $\{a_1, \ldots, a_n\}$ be a basis of $\calA$ such that $a_1$ is regular. Then
$$\calJ=\{J_{i,j}= a_i \l^j + \sigma(a_i)\l^{-j}\n 1\leqslant i\leqslant n,\, j\geqslant 1\}$$
 is a vacuum sequence in $\ti\calL_+^\tau(\calG)$.  We call the soliton hierarchy constructed from this splitting and vacuum sequence the {\it $U$-hierarchy twisted by $\sigma$\/}. 

 The space defined by \eqref{ak} is 
$$Y=\{ba_1 b^{-1}\l + v+\sigma(ba_1 b^{-1})\l^{-1}\n  b\in B,\,  v\in\calK\}.$$
The flows in the $U$-hierarchy  twisted by $\sigma$ are equations for 
$$P= (ba_1b^{-1})\l + v + \sigma(ba_1b^{-1}) \l^{-1}, \quad b:\R^2\to K, v:\R^2\to \calK.$$ 
To write down the flow equation, we need to find $$Q_j(P)=\sum_{i\leq 1} Q_{j,i}(P)\l^i$$ satisfying \eqref{bt}. 
Note that 
$$b^{-1}(\p_x+ P)b =\p_x + a_1\l + (b^{-1}v b + b^{-1}b_x) + b^{-1}\sigma(ba_1b^{-1})b \l^{-1}.$$
So we can use \eqref{bt} to solve $Q_{i, j}$ the same way as for the $U$-hierarchy (cf. \cite{TerUhl98}) and conclude that $Q_{i, j}$ are polynomial differential operators in $v, b$.   
Hence the flows in this hierarchy are PDEs.  

\ss\ni {\bf The $\frac{U}{K_1}$-hierarchy twisted by $\sigma_2$}

Let  $\frac{U}{K_1}, \frac{U}{K_2}$ denote the symmetric spaces defined by involutions $\tau, \sigma_1$ and $\tau, \sigma_2$ respectively, and $\calU_i= \calK_i + \calP_i$ the corresponding Cartan decompositions for $i=1,2$.  Assume that $\sigma_1\sigma_2= \sigma_2 \sigma_1$,
\beq\label{ah}
K_1\cap K_2= S_1\times S_2, \quad K_1= S_1\times K_1', \quad K_2= K_2'\times S_2
\eeq
as direct sum of subgroups, and there is a maximal abelian subalgebra $\calA$ in $\calP_1$ such that $\sigma_2(\calA) = \calA$.
Below we associate to such pair of symmetric spaces a splitting and a vacuum sequence

\ss
 First note that
$$K_1'\cap K_2'=\{e\}, \quad \frac{K_1}{K_1\cap K_2}= \frac{K_1'}{S_2}, \quad \frac{K_2}{K_1\cap K_2}= \frac{K_2'}{S_1},$$
and at the Lie algebra level we have
\ben
\item $\calK_1\cap \calK_2= \calS_1\oplus\calS_2$ as direct sum of Lie subalgebras,
\item 
$\calK_1'= \calS_2\oplus \calQ_1$ and 
$\calK_2':=\calS_1\oplus \calQ_2$ are the Cartan decomposition of symmetric spaces $\frac{K_1'}{S_2}$ and $\frac{K_2'}{S_1}$ respectively,
where $\calQ_1= \calK_1\cap\calP_2$ and $\calQ_2= \calK_2\cap \calP_1$,
\item  $\calK_1= \calS_1 \oplus\calK_1'=\calS_1\oplus \calS_2\oplus \calQ_1$ and $\calK_2=\calS_1\oplus\calS_2\oplus\calQ_2$ as direct sum of linear subspaces.
\een

Let $\ti L^\tau(G)$ be the group of smooth maps $f$ from the circle $S^1_{\e^{-1}}$ to $G$ satisfying $\tau(f(\bar\l))= f(\l)$.  Let  $\ti L^{\tau, \sigma_1}(G)$ be the subgroup of $f\in \ti L^\tau(G)$ fixed by the involution $\ti \sigma_1$ on $\ti L^\tau(G)$ defined by 
$$\ti \sigma_1(f)(\l)= \sigma_1(f(-\l)),$$ 
$\ti L_+^{\tau, \sigma_1}$ the subgroup of $f\in \ti L^{\tau,\sigma_1}(G)$ such that there exists a holomorphic map $h$ on $\e< |\l|< \e^{-1}$ satisfying the following conditions:
\ben
\item $h$ has boundary value $f$ on $S_{\e^{-1}}^1$ and $\hat f(\l)= \sigma_2(f(\l^{-1})$ on $S_\e^1$, 
\item $h(1)\in K_2'$.
\een
It follows from the Liouville Theorem and $K_1'\cap K_2'=\{e\}$ that $\ti L_+^{\tau, \sigma_1}(G)$ intersects $\ti L_-^{\tau, \sigma_1}(G)$ at the identity. 

 Let $\ti L_-^{\tau,\sigma_1}(G)$ denote the subgroup of $f\in \ti L^{\tau, \sigma_1}(G)$ such that there exists a holomorphic map on $\infty \geq |\l|>\e^{-1}$ with boundary value $f$ on $S^1_{\e^{-1}}$ and $f(\infty)\in K_1'$.   

The Lie algebras written in Laurent series are:
\begin{align*}
\ti\calL^{\tau, \sigma_1}(\calG)&=\left\{\xi(\l)=\sum_j \xi_j \l^j\,\big|\, \tau(\xi(\bar\l))=\xi(\l),\, \sigma_1(\xi(-\l))= \xi(\l)\right\},\\
&=\left\{ \sum_{j} \xi_j \l^j\,\big|\, \xi_j\in \calK_1 \,\, {\rm if\, } j \, {\rm is \, even,\/}\,\, \xi \in \calP_1 \, \,{\rm if\,} j\, {\rm is \, odd}\right\},\\
\ti\calL_+^{\tau,\sigma_1}(\calG)&= \left\{\sum_j \xi_j \l^j\in \ti\calL^{\tau, \sigma_1}(\calG)\,\big|\, \xi_{-j}= \sigma_2(\xi_j), \, \xi(1)\in \calK_2'\right\},\\
\ti \calL_-^{\tau, \sigma_1}(\calG)&=\left\{ \xi(\l)=\sum_{j\leq 0}\xi_j \l^j \in \ti\calL^{\tau, \sigma_1}(\calG)\,\big| \, \xi_0\in \calK_1'\right\}.
\end{align*}

\bthm
$(\ti \calL^{\tau,\sigma_1}_+, \ti \calL^{\tau,\sigma_1}_-)$ is a splitting of $\ti \calL^{\tau,\sigma}(G)$.
\ethm

\begin{proof}
Given $A=\sum_i A_i \l^i\in \ti \calL= \ti \calL^{\tau, \sigma_1}(\calG)$, we want to prove that there exist uniuqe $\xi=\sum_i \xi_i \l^i$ in $\calL_+= \ti\calL^{\tau,\sigma_1}_+(\calG)$ and $\eta=\sum_i \eta_i \l^i$ in $\calL_-=\ti \calL^{\tau,\sigma_1}_-(\calG)$ such that $A= \xi+ \eta$.  It is easy to see that 
$$
\xi_j= \bca A_j, & j >0,\\ \sigma(A_{-j}), & j<0,\eca   \qquad \eta_j= A_j- \sigma_2(A_{-j}) \, \, {\rm for\, } j<0.$$
So it remains to prove that $\xi_0$ and $\eta_0$ can be solved from $A$ uniquely.  To do this, we note that $\xi=\xi_0+ \sum_{j>0} (A_j \l^j + \sigma_2(A_j)\l^{-j})$ and $\xi_0\in \calK_1$.   Since $A\in \ti\calL$, we have  $A_j\in \calK_1$ for even $j$ and $A_j\in \calP_1$ for odd $j$.  By definition of $\calL_+$,
\beq\label{by}
\xi(1)=\xi_0+\sum_{j>0} A_j+\sigma_2(A_j)\, \in\,  \calK_2'= \calS_1 + \calQ_2.
\eeq
But
$\sum_{j>0} A_j+\sigma_2(A_j)\in \calK_2$ and $\calK_2'\subset \calK_2$ imply that $\xi_0\in \calK_2$. So $\xi_0\in \calK_1\cap \calK_2=\calS_1+\calS_2$.  Note also that 
\begin{align*}
&\sum_{j>0,\, j\, {\rm even}}A_j+ \sigma_2(A_j)\in \calK_2\cap \calK_1= \calS_1+\calS_2,\\
&\sum_{j>0,\, j\, {\rm odd}} A_j +\sigma_2(A_j)\in \calK_2\cap \calP_1= \calQ_2.
\end{align*}
Thus by \eqref{by} we have 
$$\pi_{\calS_2}(\xi_0)=- \pi_{\calS_2}\left(\sum_{j>0,\,\, j\, {\rm even}}(A_j + \sigma_2(A_j))\right).$$
Use the facts that $A_0\in \calK_1$, $A_0=\xi_0 + \eta_0$, $\xi_0\in \calS_1\oplus \calS_2$, $\eta_0\in \calK_1'= \calS_2\oplus \calQ_1$, and  $\calK_1= \calS_1 \oplus \calS_2 + \calQ_1$ to conclude that
\begin{align*}
&\pi_{\calS_1}(\xi_0)=\pi_{\calS_1}(A_0),\\
&\pi_{\calS_2}(\xi_0)=- \pi_{\calS_2}\left(\sum_{j>0,\,\, j\, {\rm even}}(A_j + \sigma_2(A_j))\right),\\
&\pi_{\calS_2}(\eta_0)= \pi_{\calS_2}(A_0-\xi_0) = \pi_{\calS_2}\left(A_0 + \left(\sum_{j>0,\,\, j\, {\rm even}}(A_j + \sigma_2(A_j))\right)\right),\\
& \pi_{\calQ_1}(\eta_0)= \pi_{\calQ_1}(A_0).
\end{align*}
\end{proof}

Let $\calA$ be a maximal abelian subalgebra in $\calP_1$, and $\{a_1, \ldots, a_n\}$ a basis of $\calA$ such that $a_1$ is regular with respect to the $\Ad(K_1)$-action.  Since $\sigma_2(\calA)\subset \calA$,
$$\calJ=\{J_{i, j}= a_i \l^{j} + \sigma_2(a_i) \l^{-j}\n 1\leq i\leq n, \, j>0 \,\, {\rm odd\/}\}$$
is a vacuum sequence. 

A computation shows that the space $Y$ defined by \eqref{ak} is 
$$Y=\{ha_1h^{-1}\l + v + \sigma_2(ha_1 h^{-1})\l^{-1}\n  h\in K_1', \, v\in \calS_1\}.$$
Use the same argument as for the $U$-hierarchy twisted by $\sigma$ to see that the flows in this hierarchy are PDEs for $h:\R^2\to K_1'$ and $v:\R^2\to \calS_1$.  

\beg\label{br} {\bf The $\frac{O(n,n)}{O(n)\times O(n)}$-hierarchy twisted by $\sigma_2$}\par
\eeg

This is the example that motivates the construction of twisted $\frac{U}{K}$- hierarchy. Let $\calG= o(n, n, \C)$, $\tau(x)= \bar x$, $\sigma_1(x)= \I_{n,n} x \I_{n,n}^{-1}$, and $\sigma_2(x)= \I_{n+1, n-1} x \I_{n+1, n-1}^{-1}$, where $\I_{m,n_0-m}$ is the $n_0\times n_0$ diagonal matrix with $ii$-th entry equals to $1$ if $i\leqslant m$ and equals to  $-1$ if $m< i\leqslant n_0$.  Then 
\begin{align*}
&\calU= o(n,n), \quad \calK_1= o(n)\times o(n), \quad \calK_2= o(n,1)\times o(n-1),\\
& \calK_1\cap \calK_2= \calS_1+\calS_2, \, {\rm where }\quad \calS_1=o(n)\times 0, \, \calS_2= 0\times o(n-1),\\
&\calK_2'= \calS_1 +(\calP_1\cap \calK_2)= o(n,1),\quad \calK_1'= 0\times o(n).
\end{align*}
So $K_1, K_2$ satisfy \eqref{ah}.  
The space 
$$\calA=\left\{\bpm 0& D\\ D&0\epm\, \bigg|\, D\in gl(n,\R)\, {\rm is\, diagonal\/}\right\}$$
is a maximal abelian subalgebra in $\calP_1$ and $\sigma_2(\calA)\subset \calA$.  Choose a basis $\{a_1, \ldots, a_n\}$ of $\calA$ such that $a_1$ is regular. Then 
$$\calJ=\{ J_{i,j}= a_i\l^j + \sigma_2(a_i)\l^{-j}\n 1\leq i\leq n, \, j\,\, {\rm odd}\}$$ 
is a vacuum sequence, and we obtain the $\frac{O(n,n)}{O(n)\times O(n)}$-hierarchy twisted by $\sigma_2$.  The flows in this hierarchy are PDE for maps $g:\R^2\to O(n)$ and $v:\R^2\to o(n)$.

\section{The first order PDE systems}

Because all flows in the soliton hierarchy associated to a symmetric spaces commute, we can put the first $n$ flows together (here $n$ is the rank of the symmetric space) to construct a first order PDE system.  Many of these systems occur in submanifold geometry.

\ss\ni {\bf The $U$-system}  {\rm (\cite{Ter97})}

We use the same notations as for the $U$-hierarchy, and put  all the flows in the $U$-hierarchy generated by $a_1\l, \ldots, a_n\l$ together. In fact, if $u(t_1, \ldots, t_n)$ satisfies \eqref{aj} for each $1\leqslant i\leqslant n$, then $v: \R^n\to \calA^\perp$ defined by $v= \ad(a_1)^{-1}(u)$ satisfies
\beq\label{ai}
[a_i, v_{x_j}]-[a_j, v_{x_i}]= [[a_i, v], [a_j, v]], \quad 1\leqslant i\not=j\leqslant n,
\eeq
where $\calA^\perp=\{x\in \calU\n (x, \calA)=0\}$ and $(\, , )$ is an ad-invariant bi-linear form on $\calU$.
System \eqref{ai} is called the {\it U-system\/} This system is independent of the choice of the basis of $\calA$ because if we change basis of $\calA$, \eqref{ai} just amounts to a linear change of coordinates of $x_1, \ldots, x_n$. So to write down the $U$-system, we can use any basis for $\calA$.  Moreover, the $U$-system has a Lax pair:
$$\o_\l= \sum_{i=1}^n (a_i\l + [a_i, v])\, dx_i,$$
 i.e., $v$ is a solution of the $U$-system if and only if $\o_\l$ is a flat connection $1$-form on $\R^n$ for all $\l\in \C$.

\ss

\ni {\bf The $\frac{U}{K}$-system}  {\rm (\cite{Ter97})}\par

\ss
We use the same notations as for the $\frac{U}{K}$-hierarchy, and $\{a_1, \ldots, a_n\}$ is a basis of a maximal abelian subalgebra $\calA$ in $\calP$.  The $\frac{U}{K}$-system is the system \eqref{ai} for maps $v:\R^n\to\calA^\perp\cap \calP$.  Its Lax pair has the same form as the $U$-system.

\ss
\ni {\bf The $U$-system twisted by $\sigma$}

\ss
We use the same notation as for the $U$-hierarchy twisted by $\sigma$.   Let  $\{a_1, \ldots, a_n\}$ be a base of $\calA$ such that $a_1$ is regular. 
The  {\it $U$-system twisted by $\sigma$\/} is the collection of flows in the $U$-hierarchy twisted by $\sigma$ generated by $J_{1, i}= a_i\l + \sigma(a_i)\l^{-1}$ for $1\leq i\leq n$, i.e., it is the equation for $g:\R^n\to B$ and $v:\R^n\to \calK$ with a Lax pair of the form
\beq\label{aq}
\o_\l= \sum_{i=1}^n ((ga_i g^{-1}) \l + v_i + \sigma(ga_i g^{-1})\l^{-1})\, dx_i,
\eeq
with $g:\R^n\to B$ and $v_i:\R^n\to \calK$.

\ss
\ni {\bf The $\frac{U}{K_1}$-system twisted by $\sigma_2$}

We use the same notations as for the $\frac{U}{K_1}$-hierarchy twisted by $\sigma$, and let $\{a_1, \ldots, a_n\}$ be a base of $\calA$.  The $\frac{U}{K_1}$-system twisted by $\sigma_2$ is the collection of commuting flows in the $\frac{U}{K_1}$-hierarchy generated by $J_{1, i}= a_i\l + \sigma(a_i)\l^{-1}$ for $1\leq i\leq n$, and this system has a Lax pair of the form
$$\o_\l= \l gDg^{-1} + v + \l^{-1} \sigma_2(gDg^{-1}),$$
where $g:\R^n\to K_1'$, $v=\sum_{i=1}^n v_i \rd x_i$ with $v_i: \R^n \to \calS_1$ for $1\leq i\leq n$, and $D = \sum_{i=1}^n a_i \rd x_i$.  

The gauge transform of $\o_\l$ by $g^{-1}$
$$g^{-1}\ast \o_\l= D\l + (v+ g^{-1} \rd g) + \l^{-1} g^{-1}\sigma_2(gDg^{-1})g$$
is flat for all non-zero $\l$. If $\frac{U}{K_1}$ has maximal rank, i.e., the rank of $\frac{U}{K_1}$ is equal to the rank of $U$, then the flatness of  $g^{-1}\ast\o_\l$ implies that there is $\xi:\R^n\to \calA^\perp\cap\calP_1$ such that 
$v_i + g^{-1} g_{x_i} =[ a_i, \xi]$. So we have 

\bprop\label{bb} 
If $\frac{U}{K_1}$ has maximal rank, then the $\frac{U}{K_1}$-system twisted by $\sigma_2$ is the equation for $(g, \xi): \R^n\to K_1' \times (\calA^\perp\cap\calP_1)$ such that 
$$\o_\l=\sum_{i=1}^n ( ga_ig^{-1} \l + \pi_{\calS_1}([a_i, \xi]) + \sigma_2(gDg^{-1})\l^{-1}) \rd x_i$$
is flat for all $\l\in \C\setminus \{0\}$, where $\pi_{\calS_1}$ is the projeciton of $\calK_1$ onto $\calS_1$ along $\calK_1'$ and $g^{-1}\rd g= \sum_{i=1}^n \pi_{\calK_1'}([a_i, \xi]) \rd x_i$.    
\eprop

\ms
\beg\label{bq} {\bf A twisted $\frac{O(n,n)}{O(n)\times O(n)}$-system}

\eeg

We use the same notations as the $\frac{O(n,n)}{O(n)\times O(n)}$ system twisted by $\sigma_2$, and choose the following basis of $\calA$:
$$a_i=\frac{1}{2}\bpm 0 & e_{ii}\\ e_{ii} &0\epm, \quad 1\leqslant i\leqslant n,$$
 By Proposition \ref{bb}, the $\frac{O(n,n)}{O(n)\times O(n)}$-system twisted by $\sigma_2$ is the system for $(A, F):\R^n\to O(n)\times gl_\ast(n)$ such that 
 \beq\label{as}
 \o_\l =\frac{\l}{2} \bpm 0& \d A^t\\ A\d & 0\epm + \bpm \d F- F^t\d &0\\ 0&0\epm +\frac{\l^{-1}}{2} \bpm 0& \d A^t J\\ JA\d& 0\epm\,
 \eeq
is flat for all $\l\in \C\setminus \{0\}$, where 
$\d=\diag(dx_1, \ldots, dx_n)$, $ J=\I_{1, n-1}$, and $gl_\ast(n)=\{(y_{ij})\in gl(n)\n y_{ii}=0, 1\leq i\leq n\}$.

\ms
\ni{\bf Construction of local solutions from $f\in \ti L_-^{\tau, \sigma}(G)$}

Since flows the $\frac{U}{K_1}$-system twisted by $\sigma_2$ commute, the formal inverse scattering for the $\frac{U}{K_1}$-system twisted by $\sigma_2$ works the same way as for flows of soliton hierarchy: 
Set $V(x)= \exp(\sum_{j=1}^n (a_j \l+ \sigma_2(a_j)\l^{-1}) x_j)$.  Given $f\in \ti L_-^{\tau, \sigma}(G)$, by Theorem \ref{bv}, there is an open subset $\calO$ of the origin in $\R^n$ such that 
$fV(x)= E(x) M(x)$ with $E(x)\in L^\tau_+(G)$ and $M(x)\in L^\tau_-(G)$ for all $x\in \calO$. Then $E$ is the frame of a solution of the the $\frac{U}{K_1}$-system twisted by $\sigma_2$.

\ms
\ni{\bf Higher flows for the $\frac{U}{K_1}$-system twisted by $\sigma_2$}

Let $a=\sum_{i=1}^n c_i a_i\in \calA$ be a regular element, and $\frac{\p}{\p x}= \sum_{i=1}^n c_i\frac{\p}{\p x_i}$. 
Let $g:\R\to K_1'$ and $\xi:\R\to \calA^\perp\cap \calP_1$ satisfying $g^{-1}g_x= \pi_{\calK_1'}([a,\xi])$, and
$P= gag^{-1}\l + \pi_{\calS_1}([a, \xi]) + \sigma_2(gag^{-1})\l^{-1}$. Let $Q_i(P)= \sum_{j\leq 1} Q_{i, j}(P) \l^j$ be the solution for \eqref{bt}.  Then
the flow equation generated by $J_{i, j}$,
$$\frac{\p P}{\p t_{i,j}}= \left[\p_x + P , \pi_+(Q_{i, j}(P))\right],$$
leaves the space of solutions of the $\frac{U}{K_1}$-system twisted by $\sigma_2$ invariant. 

\ms
\ni{\bf Symmetries of the $\frac{U}{K_1}$-system twisted by $\sigma_2$}

\ss
The group $\ti L^{\tau, \sigma}_-(G)$ acts on the space of local solutions of the $\frac{U}{K_1}$-system twisted by $\sigma_2$ by dressing action as follows:  Let $E$ be the frame of  a solution $P$ of the $\frac{U}{K_1}$-system twisted by $\sigma_2$.
Given $f\in\ti  L^{\tau, \sigma}_-(G)$, by Theorem \ref{bv}, there is an open subset $\calO$ of $\R^n$ such that we can factor $fE(x)=\ti E(x) \ti f(x)$ with $\ti E(x)\in L^\tau_+(G)$ and $\ti f(x)\in L^\tau_-(G)$ for all $x\in \calO$. 
Then $\ti E$ is the frame for a new solution $\ti P$. Moreover, $f\ast P= \ti P$ 
 defines an action of $\ti L_-^{\tau, \sigma}(G)$ on the space of local solutions of the $\frac{U}{K_1}$-system twisted by $\sigma_2$.
 If $f\in L^\tau_-(G)$ is a rational map, then the action $f\ast v$ can be computed by an algebraic formula in terms of $E$ and poles and residues of $f$.

\bs

\ni \ss{\bf The GSGE}

Let $M^n$ be a simply connected submanifold of $\R^{2n-1}$ with constant sectional curvature $-1$. Then the normal bundle $\nu(M)$ is flat and  there exist coordinates $(x_1, \ldots, x_n)$, an $O(n)$-valued map $A=(a_{ij})$, and parallel normal frames $e_{n+1}, \ldots, e_{2n-1}$ such that the first and second fundamental forms are of the form
$$\I= \sum_{i=1}^n a_{1i}^2 dx_i^2, \quad \II= \sum_{i=1, j=2}^n a_{1i} a_{ji} dx_i^2 e_{n+j-1}.$$

We use the method of moving frame to write down the Gauss-Codazzi equation for these immersions:  Set 
$$F=(f_{ij}), \quad  f_{ij}= \bca\frac{(a_{1i})_{x_j}}{a_{1j}}, & i\not= j,\\ 0,& i=j,\eca,$$ and $w_{ij}= f_{ij}dx_i - f_{ji} dx_j$ for $1\leq i, j\leq n$. The Gauss-Codazzi equation is
\beq\label{bp}
\bca dw+ w\wedge w= -\frac{1}{2} \d A^t e_{11}A\d, & {\rm where\, } w= \d F- F^t \d,\\
A^{-1}dA= \d F^t - F\d.\eca
\eeq
This is the GSGE (\cite{TenTer80, Ter80}).  It has a Lax pair (\cite{AbBeTe86}):
\beq\label{bs}
\W_\l= \bpm \d F- F^t\d & \frac{\l}{2} \d A^t + \frac{\l^{-1}}{2} \d A^t J\\ \frac{\l}{2} A\d + \frac{\l^{-1}}{2} JA\d& 0\epm.
\eeq
Note that  \eqref{bs} is the same as \eqref{as}. So we have

\bthm
The GSGE is the twisted $\frac{O(n,n)}{O(n)\times O(n)}$-system given in Example \ref{bq}.
\ethm

Since the twisted $\onn$-system comes from a splitting of a subalgebra of a loop algebra, we obtain  higher flows and symmetries of  GSGE.

\end{document}